\documentclass[12pt]{amsart}
\usepackage{amsmath}
\usepackage{amsthm}
\usepackage{amscd}

\newcommand{\Q}{{\mathbb Q}}
\newcommand{\R}{{\mathbb R}}
\newcommand{\Z}{{\mathbb Z}}
\newcommand{\C}{{\mathbb C}}

\newcommand{\A}{{\mathbf A}}

\newtheorem{theorem}{Theorem}
\newtheorem{corollary}{Corollary}
\newtheorem{proposition}[theorem]{Proposition}

\theoremstyle{definition}
\newtheorem{conjecture}{Conjecture}

\newtheorem{definition}{Definition}
\numberwithin{definition}{section}

\numberwithin{equation}{section}

\theoremstyle{remark}
\newtheorem*{remark}{Remark}
\newtheorem*{acknowledgement}{Acknowledgement}

\begin{document}

\title[On isospectral arithmetical spaces]{On 
  isospectral arithmetical spaces}

\author{C.~S.~Rajan}

\address{Tata Institute of Fundamental  Research, Homi Bhabha Road,
Bombay - 400 005, INDIA.}  \email{rajan@math.tifr.res.in}


\begin{abstract}
We study the relationship between the arithmetic and the spectrum of
the Laplacian for manifolds  arising from congruent arithmetic
subgroups of $SL(1,D)$, where $D$ is an indefinite quaternion division
algebra defined over a number field $F$.
 
We give new examples of isospectral but non-isometric compact,
arithmetically defined varieties, generalizing the class of examples
constructed by Vigneras. These examples are based  on an interplay
between the simply connected and adjoint group and  depend explicitly
on the failure of strong approximation for the adjoint group. The
examples can be considered as a geometric analogue and  also  as an
application of the concept and results on  $L$-indistinguishability
for $SL(1,D)$ due to  Labesse and Langlands.

We verify that the Hasse-Weil zeta functions  are equal for the
examples of isospectral pair of Shimura varieties we construct giving
further evidence for an archimedean analogue of Tate's conjecture,
which expects  that the spectrum of the Laplacian determines the
arithmetic of such spaces.

\end{abstract}
 
\maketitle

\section{Introduction}  
Let $M$ be a compact, connected Riemannian manifold.  The spectrum of
$M$ consists of the collection of non-zero eigenvalues counted with
multiplicity of the Laplace-Beltrami operator acting on the space of
smooth functions on $M$.  The inverse spectral problem is the
investigation of the properties of the manifold that can be recovered
from a knowledge of the spectrum.   Two compact connected Riemannian
manifolds $M_1$ and $M_2$ are said to be isospectral if the spectrums
of $M_1$ and $M_2$ coincide.  The first question that arises in
inverse spectral theory is whether the spectrum determines a
Riemannian manifold upto isometry.   This is known to be false and the
first counter-examples were constructed by  Milnor in the case of flat
tori. When the spaces are compact hyperbolic surfaces, counterexamples
were constructed by  by Vigneras \cite{V}  and later by Sunada
\cite{S}.  In this paper, our aim is to investigate the inverse
spectrum problem in the arithmetical context  of compact, locally
symmetric manifolds, especially those arising from co-compact
congruence lattices in $G=SL(1,D)$, where $D$ is an indefinite
quaternion division algebra over a number field.

The main observation is that of introducing the conjugation of a
compact open subgroup $K$ of $G(\A_f)$  by an element of the adelic
points of the adjoint group, and to observe that this process yields
spectrally indistinguishable manifolds.  We establish this observation
(see Thereom \ref{adelicconjtheorem}) under some hypothesis on the
nature of $D$ and the compact open subgroup $K$. This result can be
considered as a geometric analogue of the study initiated by Langlands
together with Labesse and Shelstad on the problems arising out of the
difference between stable conjugacy and conjugacy: in \cite{LL}, the
multiplicity of the representations $\pi^g$ in the space of cusp forms
on $SL(2, \A_F)$ as a function of $g$ is studied,  where $\pi$ is  a
representation of $SL(2,\A_F)$ and $g\in GL(2, \A_F)$.

As a corollary to the above observation, we give examples of
isospectral but non-isometric compact, locally symmetric  manifolds as
above, generalizing  the class of isospectral but non-isometric
hyperbolic surfaces and three dimensional manifolds  constructed by
Vigneras \cite{V}.  For $K$ as above, the condition for the existence
of  non-isometric spaces is that the finite adelic points of the
adjoint group is not exhausted by the union of translates of the
normalizer of the image of  $K$ by rational elements. This always
happens if the normalizer of $K$ is sufficiently small by the  failure
of strong approximation in the adjoint group.   In contrast to the
method of this paper, Vigneras' method is  a combination of geometry
and arithmetic. It depends  on the relationship of the spectrum of the
Laplacian with the length spectrum for compact hyperbolic
surfaces. They work when the lattices are associated to maximal
orders,  and is difficult to generalize to more general groups or to
more general lattices.
 
One of the motivating questions for us has been the heuristic that the
  spectrum of invariant operators should determine the  arithmetic of
  compact locally symmetric manifolds arising from arithmetic
  lattices.  For the family of arithmetic curves associated to
  congruent lattices coming from quaternion division algebras over
  totally real fields, Shimura constructed a theory of canonical
  models \cite{Sh}. In particular, these spaces can be defined over
  number fields. The question was raised in \cite{PR} that if two
  arithmetically defined curves are isospectral, then are the
  Hasse-Weil zeta functions of suitably defined canonical models
  equal?  The conjecture is based on the analogy that the Laplacian
  can be considered as the Frobenius at infinity. Further, the theorem
  of Faltings proving Tate's conjecture asserts that  if the
  eigenvalues of the Frobenius elements acting on the $l$-adic
  cohomology groups of two smooth, projective curves defined over a
  number field coincide, then the Jacobians of the curves are
  isogenous. Thus the conjecture can be considered as an archimedean
  analogue of Tate's conjecture.  In \cite{PR}, this conjecture was
  verified in the examples constructed by Sunada and Vigneras.  Using
  the basic fact in the theory of canonical models that adelic
  conjugation corresponds to Galois conjugation of the canonical
  model, we verify the conjecture that the spectrum determines the
  Hasse-Weil zeta function for  the examples constructed in this
  paper.  On the other hand, the results of this paper can be used to
  reformulate the conjecture for the non-connected Shimura varieties
  rather than the connected ones.

To  prove  Theorem \ref{adelicconjtheorem}, by the  generalized Sunada
 criterion for isospectrality due to DeTurck and Gordon, it is
 sufficient to  show that the lattices corresponding to $K$ and it's
 conjugate are representation equivalent. Using strong approximation
 for $SL(1,D)$,  the proof of the representation equivalence of the
 two lattices is an application of the formula for the multiplicity of
 a representation of $G(\A)$ in the space of cusp forms on $G(\A)$ due
 to Labesee and Langlands.

\begin{acknowledgement} 
I sincerely thank Cefipra for supporting my travel and stay at
  Universit\'{e} de Paris VII (Jussieu) during May of 2002, during
  which period part of this work was carried out. I thank Bruno Kahn
  and Paul Gerardin for their hospitality and a comfortable stay at
  Paris, especially to P. Gerardin for graciously providing me with
  office space during my stay in Paris.

My sincere thanks to  J. P. Labesse and A. Raghuram for useful
discussions on $L$-indistinguishability for $SL(2)$, and to D. Prasad
for many useful discussions.
\end{acknowledgement}

\section{Sunada's criterion and representation equivalence}
Let $G$ be a Lie group, and let $\Gamma$ be a cocompact lattice in
$G$. The existence of a cocompact lattice implies that $G$ is
unimodular. Let $R_{\Gamma}$ denote the right regular representation
of $G$ on the space $L^2(\Gamma\backslash G)$ of square integrable
functions with respect to the projection of the Haar measure on the
space $\Gamma\backslash G$:
\[(R_{\Gamma}(g)f)(x)=f(xg)\quad f\in  L^2(\Gamma\backslash G), ~g,
x\in G.\]

\begin{definition} Let $G$ be a Lie group and $\Gamma_1$ and
  $\Gamma_2$ be two co-compact lattices in $G$. The lattices
$\Gamma_1$ and $\Gamma_2$ are said to be {\em representation
equivalent in $G$}  if the regular representations $R_{\Gamma_1}$ and
$R_{\Gamma_2}$ of $G$ are isomorphic.
\end{definition}

We have the following generalization of Sunada's criterion for
isospectrality proved by DeTurck and Gordon in \cite{DG}:
\begin{proposition}\label{sunada}
Let $G$ be a Lie group  acting on the left as isometries of a
Riemannian manifold $M$.  Suppose $\Gamma_1$ and $\Gamma_2$ are
discrete, co-compact  subgroups of $G$ acting freely and properly
discontinuously on $M$, such that the quotients $\Gamma_1\backslash M$
and $\Gamma_2\backslash M$ are compact Riemannian manifolds. If  the
lattices  $\Gamma_1$ and $\Gamma_2$ are representation equivalent in
$G$, then  $\Gamma_1\backslash M$ and $\Gamma_2\backslash M$ are
isospectral for the Laplacian acting on the space of smooth functions.
\end{proposition} 

\begin{remark} When $G$ is a finite group, this result is due to
  Sunada \cite{S}. Sunada's construction is based on an analogous
  construction in algebraic number theory: let  $G$ be the Galois
  group of a finite Galois extension $L$ over  the rationals. Suppose
  there exists two subgroups  $\Gamma_1$ and $\Gamma_2$  of $G$  such
  that the  representations $R_{\Gamma_1}$ and $R_{\Gamma_2}$ are
  isomorphic. Then the invariant  fields $L^{\Gamma_1}$ and
  $L^{\Gamma_2}$ have the same Dedekind zeta function.

The conclusion in Proposition \ref{sunada} can be strengthened to
  imply  `strong isospectrality', i.e., the spectrums coincide for
  natural self-adjoint  differential operators besides the Laplacian
  acting on functions. In particular, this implies the isospectrality
  of the Laplacian acting on $p$-forms.

A naive expectation is that suitable strong isospectrality assumptions
on a pair of compact locally symmetric spaces should determine the
`arithmetic' associated to these spaces (although it does  not
determine the geometry of such spaces).
\end{remark}

\section{Adelic conjugation of lattices} 
Let $F$ be a number field and let $D$ be  a quaternion algebra  over
$F$.
\[ \text{Let}\quad  G=SL_1(D), \quad \tilde{G}=GL_1(D), \quad
\bar{G}=PGL_1(D).\]  Let $K_{\infty}$ be a maximal compact subgroup of
$G_{\infty}$ and let $M=G_{\infty}/K_{\infty}$ be the non-compact
symmetric space corresponding to $G$. Let $K$ be a compact open
subgroup of $G(\A_f)$.  Denote by  $\Gamma_K=G_{\infty}K\cap G(F)$ the
co-compact lattice  in $G_{\infty}=G(F\otimes \R)$ corresponding to
$K$.  If $K$ is sufficiently small, then  $\Gamma_K$ will be a
torsion-free lattice acting freely and properly discontinuously on
$M$.  The quotient space $M_K=\Gamma_K\backslash M$ will then be  a
compact, locally symmetric space.

In this section,  we assume that $D$ and $K$ satisfy the following
hypothesis:
\begin{description}
\item[H1] $D$ is a division algebra, and there is at least one  finite
place $v_0$ of $F$  at which $D$ is ramified, i.e.,  $D\otimes
F_{v_0}$ remains a division algebra.

\item[H2] $D$ is indefinite: there is at least one archimedean place
$v$ of $F$, at which  $D\otimes_F F_v\simeq M_2(\R)$,  where for a
place $v$ of $F$,  $F_v$ denotes the completion of $F$ at $v$.

\item[H3] There is a factorisation of the form,
\[ K=K_{v_0}K^{v_0},\]
where $K_{v_0}$ is a compact,  open subgroup of $G(F_{v_0})$,  and is
a normal  subgroup of $D^*_{v_0}$. The group $K^{v_0}$ has no $v_0$
component, i.e., for any element $x\in K^{v_0}$, the $v_0$ component
$x_{v_0}=1$.
\end{description}

The heart of this paper is the following theorem, whose proof is given
in Section \ref{proof}:
\begin{theorem} \label{adelicconjtheorem}
Let $D$ be a quaternion division algebra over a number field $F$ and
let $K$ be a compact open subgroup of $G(\A_f)$ satisfying the
properties {\bf H1, H2, H3} given above.  Then for any element  $x\in
GL_1(D)(\A_f)$, the lattices $\Gamma_K$  and $\Gamma_{K^x}$ are
representation equivalent in  $G_{\infty}$.
\end{theorem}
\begin{remark}
 The construction can be considered as a geometric analogue of concept
  of $L$-indistinguishability arising out of the difference between
  stable conjugacy and conjugacy  due to Labesse, Langlands and
  Shelstad \cite{LL, Shl}. Let $\pi$ be a representation of $G(\A)$
  and $g\in \tilde{G}(\A)$. The representations $\pi$ and $\pi^g$ are
  said to be $L$-indistinguishable (see \cite{LL} and Section 5).  In
  \cite{LL}, Labesse and Langlands consider the multiplicity with
  which the conjugate representations $\pi^g$ occur in the space of
  automorphic representation of $G(\A)$ as a  function of $g$. We can
  call the lattices $\Gamma_K$ and $\Gamma_{K^x}$ as `stably
  conjugate' lattices, and the theorem says that stably conjugate
  lattices gives raise to spectrally indistinguishable manifolds. But
  the terminology can be misleading, as it is not clear that any
  element of $\Gamma_K$ is stably conjugate to an element of
  $\Gamma_{K^x}$ upto an element of the center of $G$. This fact is
  true for cases of surfaces and three dimensional manifolds after
  Theorem \ref{adelicconjtheorem}, by comparing the spectrum with the
  length spectrum which is given by the absolute value of the trace
  (see \cite{V}).

 More generally, let $\Gamma_1, ~\Gamma_2$ be  two  lattices in
  $H(\R)$, where $H$ is a reductive algebraic group over $\R$. Call
  $\Gamma_1$ and $\Gamma_2$ to be weakly stably conjugate if every
  element $\gamma_1\in \Gamma_1$ is stably conjugate (i.e., conjugate
  inside the group $H(\C)$) to an element of $\Gamma_2$ upto a central
  element of $H(\R)$.  The following question can be raised: if
  $\Gamma_1$ and $\Gamma_2$ are weakly stably conjugate cocompact
  lattices in $G(\R)$, are they representation equivalent? It is
  possible that the question may have an affirmative answer with
  either a weaker conclusion: given an irreducible, unitary
  representation $\pi$  of $G(\R)$, the sum of the multiplicities of
  the representations in the (Arthur) $L$-packet of $\pi$ occuring in
  $R_{\Gamma_1}$ and $R_{\Gamma_2}$ are equal;  or with a stronger
  hypothesis: the stable conjugacy classes are counted with an
  appropriate notion of multiplicity.

But as far as producing examples of isospectral spaces are concerned,
 the problem remains of producing  examples of pairs of non-conjugate
 but weakly stably conjugate lattices.
\end{remark}

\begin{remark} 
The proof of the  theorem given here uses adelic methods and depends
quite crucially on hypothesis {\bf H3} that  $K$ is required to
satisfy.  For example, an analogous question  can be raised  for the
cuspidal spectrum of $SL(2)$ for which we do not know an answer: let
$K$ be a compact open subgroup of $SL_2(\A_f)$. For any $x\in
GL_2(\A_f)$ are the cuspidal spectrums of $L^2(\Gamma_K\backslash
SL_2(F\otimes \R))$ and $L^2(\Gamma_{K^x}\backslash SL_2(F\otimes
\R))$ equal? In this situation, hypothesis {\bf H3} cannot be assumed
as the groups $SL_2(k)$ are simple modulo center for any infinite
field $k$.
\end{remark}

\subsection{Isospectral spaces}
The following corollary gives examples of isospectral but
non-isometric compact Riemannian manifold:
\begin{corollary}\label{isospnotisom} With the notation as in Theorem
  \ref{adelicconjtheorem}, assume further that $K$ is small enough so
  that $\Gamma_K$ and $\Gamma_{K^x}$  are torsion-free.  Let
  $N(\bar{K})$ denote the normalizer of $\bar{K}$ in
  $\bar{G}(\A_f)$. Suppose $x$ is an element in $\tilde{G}(\A_f)$ such
  that $\bar{x}$ does not belong to the set $N(\bar{K})\bar{G}(F)$.
  Then $X_K$ and $X_{K^x}$ are (strongly) isospectral, but are not
  isometric.
\end{corollary} 
\begin{proof} 

Suppose on the contrary, that $X_K$ and $X_{K^x}$ are isometric. Then
there exists $\bar{g}\in \bar{G}(\R)$ such that
\[ \bar{g}^{-1}\overline{\Gamma}_{K^x}\bar{g}=\overline{\Gamma}_K,\]
where $\overline{\Gamma}_K, ~\overline{\Gamma}_{K^x}$ is the
projection of ${\Gamma_K}, ~{\Gamma_{K^x}}$ to
$\overline{G}_{\infty}$. Since the lattices ${\Gamma_K}$ and
${\Gamma_{K^x}}$ are arithmetic and commensurable, it follows  by a
theorem of Margulis that $\bar{g}\in \bar{G}(F)$. Since the kernel of
the projection map  $\tilde{G}\to G$ is a  split torus, by Hilbert
Theorem 90, there is an element $\tilde{g}\in \tilde{G}(F)$ satisfying,
\[  \tilde{g}^{-1}{\Gamma_{K^x}}\tilde{g}={\Gamma_K}.\]
Since $D$ is indefinite, the lattice $\Gamma_K$ sitting diagonally
inside $G(\A_f)$ is  dense in $K$, by the strong approximation theorem
for $SL(1,D)$. Hence
\[ \tilde{g}^{-1}K^x\tilde{g}=K.\]
Projecting to $\bar{G}$, we obtain
\[ \bar{g}^{-1}\bar{x}^{-1}\bar{K}\bar{x}\bar{g}=\bar{K}, \]
where $\bar{K}$ denotes the image of $K$ in $\bar{G}(\A_f)$.  This
implies that $\bar{x}\in N(\bar{K})\bar{G}(F)$, contradicting our
choice of $\bar{x}$.
\end{proof}

\begin{remark} The 
normalizer $N(\bar{K})$ is a compact open subgroup of $\bar{G}(\A_f)$.
By the failure of strong approximation for the adjoint group
$PGL_1(D)$, the hypothesis that $\bar{x}$ does not belong to the
double coset $N(\bar{K})\bar{G}(F)$ is satisfied provided $N(K)$ is
small enough. For example $K$ can be taken to be sufficiently deep
Hecke congruence subgroups corresponding to a finite collection of
split primes for $D$.
\end{remark}
 
\begin{remark} These examples generalize the class of examples
  constructed by Vigneras of isospectral but non-isometric compact
  Riemann surfaces. Vigneras' method  is to relate  the spectrum to
  the length spectrum which is well understood for compact hyperbolic
  surfaces and three folds. For the lattices arising from maximal
  orders in $D$, it is possible to obtain explicit formula for the
  multiplicities of lengths of periodic geodesics. Using these
  formulae, the above corollary follows in dimensions two and three
  and for lattices arising from maximal orders. The failure of strong
  approximation is reflected in the fact that the lattices correspond
  to non-conjugate maximal orders inside the quaternion algebra. In
  particular, it is required that the underlying field has a
  non-trivial class number. These methods seem difficult to generalize
  to higher dimensions and to general lattices not necessarily arising
  from maximal orders. In contrast, the above corollary is applicable
  even for indefinite quaternion division algebras over rationals.
\end{remark}

\section{Spectral  and Arithmetical Equivalence} 
In this section, we take $F$ to be a totally real number field and $D$
a quaternion division algebra over $F$. Let $\tau_1, \cdots, \tau_r$
be the real embeddings corresponding to archimedean places of $F$ at
which $D$ splits, and let $F'$ be the reflex field of $(F, \tau_1,
\cdots, \tau_r)$.  For any compact open subgroup $\tilde{K}$ of $
\tilde{G}(\A_f)$, let $F_{\tilde{K}}$ be the abelian extension of $F'$
defined as in \cite[page 157]{Sh}. Let  $K=\tilde{K}\cap G(\A_f)$. In
\cite{Sh}, Shimura defines a canonical model for the spaces $M_K$ over
the field $F_{\tilde{K}}$. As a consequence of the main theorem of
canonical models in \cite[Theorem 2.5, page 159, Section 2.6]{Sh} and
Theorem \ref{adelicconjtheorem},  we obtain the following corollary
providing more evidence in support of the conjectures made in
\cite{PR} that the spectral zeta function determines the arithmetical
zeta function:
\begin{corollary} Let  $F$ be  a totally real number field and 
$\tilde{K}$ be a compact open subgroup of $ \tilde{G}(\A_f)$. Assume
  that $\tilde{K}$ satisfies the hypothesis of Theorem
  \ref{adelicconjtheorem} and is such  that the lattices
  $\Gamma_{K^x}$ are torsion-free. Then the spaces $M_{K^x}$ for $x\in
  \tilde{G}(\A_f)$ are isospectral and have the same Hasse-Weil zeta
  function for the canonical model defined by Shimura.

If $D$ is ramified at all real places except one, then the Jacobians
of $M_{K}$ and $M_{K^x}$ are conjugate by an automorphism of
$\bar{\Q}$.
\end{corollary} 
 
There are a couple of remarks to be made about the conjecture made in
\cite{PR}. If we fix our attention on $K$, there can be more than one
choice of $\tilde{K}$ with $K=\tilde{K}\cap G(\A_f)$. Thus the choice
of a canonical model on the space $M_K$ is not uniquely determined by
$K$, and we can expect the equality of the zeta functions only after a
finite extension of the base field.

A second observation is that the zeta function of Shimura varieties
and it's relationship to that of $L$-functions of automorphic
representations are  better behaved for the non-connected Shimura
varieties associated to inner forms of  $GL(2)$, than for the
connected Shimura varieties.  Consider the space
\[ M_{\tilde{K}}= \tilde{G}(F)_+\backslash M\times \tilde{G}(\A_f)/\tilde{K}\]
 where $\tilde{G}(F)_+$ consists of those elements with totally
positive determinant. The non-connected space $M_{\tilde{K}}$ has a
canonical model in the sense of Shimura and Deligne \cite{D} over the
reflex field $F'$,   The connected components of $M_{\tilde{K}}$ are
of the form $M_{K^x}$ as $x$ ranges over $\tilde{G}(\A_f)$. We make
the following conjecture:
\begin{conjecture} Let $\tilde{K}_1, ~\tilde{K}_2$ be compact open
  subgroups of $\tilde{G}(F)$, such that the lattices $\Gamma_{K_1}$
  and $\Gamma_{K_2}$ are torsion-free modulo their centers. If the
  spaces $M_{K_1}$ and $M_{K_2}$ are isospectral, then the Hasse-Weil
  zeta functions of the canonical models of  $M_{\tilde{K}_1}$ and
  $M_{\tilde{K}_2}$ are equal.
\end{conjecture}
The conjecture is well posed thanks to Theorem
\ref{adelicconjtheorem}. This clarifies and makes precise the
arithmetical conjecture made in \cite{PR}.

\begin{remark} The initial hope behind the conjectures made in
  \cite{PR}, was that isospectral compact Riemann surfaces (with
  respect to the hyperbolic metric) have isogenous Jacobians. This
  naive conjecture proves to be false for the example of Shimura
  curves constructed above, as the respective Jacobians become
  isogenous after twisting by an automorphism of $\bar{\Q}$.

It is tempting to conjecture that this will be the exception:  if two
  compact Riemann surfaces are isospectral, then the Jacobian of one
  is isogenous to a conjugate of the Jacobian of the other by an
  automorphism $\sigma\in {\rm Aut}(\C/\Q)$, where $\sigma$ preserves
  the spectrum of the Riemann surface. If morever $\sigma$ is not
  identity or of order two, then the pair of Riemann surfaces  arise
  from an arithmetical context, i.e., are Shimura curves as considered
  in this paper.

\end{remark}

\section{Multiplicity formula  of Labesse and Langlands}
Our aim in this section is to recall the multiplicity formula of
Labesse and Langlands \cite{LL}, giving the multiplicity $m(\pi)$ with
which a representation $\pi$ of $SL_1(D)(\A_F)$ occurs in the space of
cusp forms on $SL_1(D)(\A_F)$.  This formula will be  required for the
proof of Theorem \ref{adelicconjtheorem}.  The references for the
material contained in this section are \cite{LL} and \cite{Shl}.

\subsection{Langlands parameters}
Let $F$ be a global or a local field of characteristic zero, $W_F$ be
the Weil group of $F$.  and let $W'_F$ be the Weil-Deligne group of
$F$:
\[ W_F'=\begin{cases}  W_F\times SU(2)\quad \text{if $F$
is non-archimedean}\\  W_F\quad \text{if $F$ is archimedean.}
\end{cases}\] 
By abelian reciprocity, we will identify $W_F^{ab}$ with $C_F$, where
\[C_F=\begin{cases} F^* \quad \text{if $F$ is local}\\
J_F/F^* \quad \text{if $F$ is global},
\end{cases}
\]
where for a global field $F$, $J_F$ is the idele group associated to
$F$.

For a connected  reductive group $H$ defined over $F$, let $^LH^0$ be
the connected reductive group over $\C$ whose root datum is the dual
root datum of the root datum of $H$ over $\bar{F}$.  Let $^LH$ denote
the Langlands dual group defined over the complex numbers associated
to $H$ over $F$. The Langlands dual group $^LH$ is defined as a
semidirect product,
\[ ^LH={^LH}^0\times W_F'.\]
Let $\Phi(H)$ be the set of equivalence classes of admissible
homomorphisms of $W_F'$ into $^LH$.

In our case, where $H$ is an inner form of $SL_2$, the connected
component ${^LH}^0$ is isomorphic to $PGL_2(\C)$. Further, the
Langlands dual group for $GL(2)$ defined over $F$ is just the direct
product $GL(2,\C)\times W'_F$. Let $p:{^L\tilde{G}}\to {^LG}$ denote
the natural projection map.  By a lemma of Langlands (or more
generally by the results of \cite{Lab}, \cite{Ra}), any admissible
homomorphism $\phi\in \Phi(G)$ admits a lifting $\tilde{\phi}:W_F'\to
^L\tilde{G}$, such that $p\circ \tilde{\phi}=\phi$.

The particular class of parameters that are of interest to us are
those for which $\phi$ is an induced representation.  Let $E$ be a
quadratic extension of $F$, and  $\theta$ a character  of $W_E$.
Consider the  induced representation
\[ I(\theta)=\text{Ind}_{W'_E}^{W'_F}(\theta):W_F\to {^L\tilde{G}},\]
where we consider $\theta$ as a character of   $W_E'$ trivial on
 $SU(2)$.

Denote by $\sigma$ the Galois conjugation on $E$ over $F$, and by
${\theta}^{\sigma}$ the character
${\theta}^{\sigma}(\gamma)=\theta(\sigma(\gamma))$ for  $\gamma\in
C_E$, where we consider $\theta$ as a character of $C_E\simeq
W_E^{ab}$.    The representation $I(\theta)$ is irreducible  precisely
when $\theta^{\sigma}\neq \theta$.  If $F$ is a local field and $G$ is
anisotropic over $F$, then the parameter $I(\theta)$ is admissible for
$G$ if and only if it is irreducible.

When $E/F$ is  a quadratic extension of global fields, $D$ a
quaternion division algebra over $F$, and  $\theta$  an idele class
character of $J_E/E^*$, the notion of admissibility is defined as
follows:

\begin{definition}  An idele class character $\theta$ of
$J_E/E^*$  is said to be {\em admissible} with respect to the
quaternion algebra $D$, if $\theta^{\sigma}\neq \theta$ and the
following is satisfied: let $v$ be a place of $F$ at which $D$ does
not split. Then the place $v$ should either be inert or ramified in
$E$. Let $w$ be the unique  place of $E$ dividing $v$. We further
require that the local component $\theta_w$ satisfies
$\theta_w^{\sigma}\neq \theta_w$.
\end{definition}

 Let $S_{\phi}$ be the centralizer in $^LG^0$ of the image of $\phi$,
$S_{\phi}^0$ be the connected component of the identity in $S_{\phi}$,
and let $C_{\phi}= S_{\phi}/S_{\phi}^0$ be the group of connected
components of $S_{\phi}$.  The group  $C_{\phi}$ is a finite abelian
group and is isomorphic to either $1, ~\Z/2\Z$ or $\Z/2\Z\times
\Z/2\Z$.

\subsection{Local theory} We recall the local theory. Let $F$ be
a local field of characteristic zero and $D$ be a quaternion algebra
over $F$. Let $\pi$ be a representation of $G(F)$. For $g\in
\tilde{G}(F)$, let $\pi^g$ denote the conujgate representation defined
by $\pi^g(h)=\pi(ghg^{-1})$ for $h\in G(F)$.

\begin{definition} Two irreducible, admissible representations $\pi_1$
  and $\pi_2$ of $G(F)$ are said to be {\em $L$-indistinguishable} if
  $\pi_2\simeq \pi_1^g$ for some $g\in \tilde{G}(F)$. Equivalently if
  both  $\pi_1$ and $\pi_2$ occur in the restriction of an
  irreducible, admissible  representation $\tilde{\pi}$ of
  $\tilde{G}(F)$ to $G(F)$.

Given an  irreducible, admissible  representation $\pi$ of
 $\tilde{G}(F)$ (or $\tilde{\pi}$ of $\tilde{G}(F)$), the {\em
 $L$-packet $L(\pi)$ of $\pi$} is defined to be the collection of all
 irreducible, admissible  representation of $\tilde{G}(F)$ which are
 $L$-indistinguishable from $\pi$.
\end{definition}

Since the index of $F^*G(F)$ inside $\tilde{G}(F)$ is finite,  any
irreducible, admissible  representation $\tilde{\pi}$ of
$\tilde{G}(F)$ decomposes as a finite direct sum of irreducible,
admissible representations of $G(F)$. Hence $L$-packets are of finite
cardinality. It can be seen that the cardinality of an $L$-packet is
either $1, ~2~\text{or}~4$.

If $\tilde{\phi}$ is an admissible parameter for $\tilde{G}$,  let
$\tilde{\pi}(\tilde{\phi})$ be the  irreducible, admissible
representation of $\tilde{G}(F)$ associated to $\tilde{\phi}$ by the
local Langlands correspondence. Let $L(\phi)$ be the $L$-packet given
by the collection of irreducible, admissible representations occurring
in the decomposition of $\tilde{\pi}(\tilde{\phi})$ to $G(F)$. This
depends only on the parameter $\phi$.

Fix an additive character $\psi$ of $F$.  Depending on the choice of
$\psi$, Labesse and Langlands define a pairing,
\begin{equation}\label{pairing}
 <~,~>~:~ S_{\phi}/S_{\phi}^0\times L(\phi)\to \Z.
\end{equation}
We make a provisional definition:
\begin{definition} The pairing \ref{pairing} is said to be
  non-degenerate if the pairing idenitifies $L(\phi)$
(non-canonically) as the dual group of $S_{\phi}/S_{\phi}^0$, i.e.,
$L(\phi)$ and $C_{\phi}$ are of equal cardinality and the map sending
an element $\pi\in L(\phi)$ to the character $s\mapsto <s,\pi>, ~s\in
C_{\phi}$ is a bijection from $L(\phi)$ to the dual group of
characters $\hat{C}_{\phi}$ of $C_{\phi}$.
\end{definition} 
We have,
\begin{proposition} \label{pairing:nondeg}
Let $F$ be a non-archimedean local field of characteristic zero. The
pairing \ref{pairing} is non-degenerate, except when $D$ is a
quaternion division algebra over $F$, the parameter $\tilde{\phi}$ is
induced of the form $I(\theta)$, where $\theta$ is a character of a
quadratic extension $E$ of $F$ such that $\theta^{\sigma}/\theta$ is a
non-trivial quadratic character.
\end{proposition}
In the exceptional case,  the $L$-packet consists of one
 representation $\pi$ occurring with multiplicity two in the
 restriction of $\tilde{\pi}(I(\theta))$ to $G(F)$, and
\[C_{\phi}\simeq \Z/2\Z\times \Z/2\Z.\]
The pairing is given by,
\begin{equation}\label{pairing:exceptional}
<1,\pi>=2, \quad   <s,\pi>=0, \quad s\neq 1, \quad s\in C_{\phi}.
\end{equation}

For the proof of Theorem \ref{adelicconjtheorem}, the following case
of the pairing \ref{pairing} is crucial: let  $D$ be  a quaternion
division algebra over $F$, $E$ a quadratic extension of $F$, and
$\theta$ a character of $E^*$ such that $\theta/{\theta}^{\sigma}$ is
non-trivial and not quadratic. Then the $L$-packet
$L(\phi)=\{\pi(I(\theta))^+, \pi(I(\theta))^-\}$ consists of two
representations and $C_{\phi}\simeq \{1,\epsilon\}$ is of cardinality
two, such that the pairing is given by
\begin{equation}\label{pairing:crucial}
\begin{array} {ll}
 <1,\pi(I(\theta))^{+}>=1, &\quad   <1,\pi(I(\theta))^{-}> =1\\
 <\epsilon,\pi(I(\theta))^+> =1, & \quad <\epsilon,\pi(I(\theta))^->
 =-1.
\end{array}
\end{equation}

\subsection{Global theory}
Let $F$ be a number field, and $D$ be a quaternion division algebra
over $F$.  If $\Pi$ is any irreducible, admissible representation of
$G(\A_F)$,  then $\Pi$ can be decomposed
\[ \Pi=`\otimes'_{v\in \Sigma_F}\Pi_v,\]   
as a restricted tensor product of the local components $\Pi_v$ of
$\Pi$ at the places $v$ of $F$. The local component $\Pi_v$ is an
irreducible, admissible representation of $G(F_v)$, such that at
almost all finite places $v$ of $F$, the representation $\Pi_v$ is an
unramified representation of $G(F_v)$.

\begin{definition}
Two irreducible, admissible representations $\Pi_1$ and $\Pi_2$ of
  $G(\A_F)$ are said to be {\em $L$-indistinguishable} if for all
  places $v$ of $F$, the local components  $\Pi_{1,v}$ and $\Pi_{2,v}$
  are $L$-indistinguishable, and $\Pi_{1,v}\simeq  \Pi_{2,v}$ at
  almost all places of $F$.

 The {\em $L$-packet $L(\Pi)$} of an  irreducible, admissible
 representation $\Pi$ of $G(\A_F)$ is  the collection of all
 irreducible, admissible  representation of $G(\A_F)$ which are
 $L$-indistinguishable from $\Pi$.  Equivalently it is the orbit of
 $\Pi$ under the action of $\tilde{G}(\A_F)$.
\end{definition}

\subsubsection{Dihedral representations}
Let $E/F$ be a quadratic extension of $F$, and $\theta$ a character of
$C_E$ satisfying $\theta^{\sigma}\neq \theta$. Generalizing earlier
construction of Hecke and Mass, Jacquet and Langlands associate a
cuspidal automorphic representation $\tilde{\Pi}(I(\theta))$ of
$GL_2(\A)$ such that at any place $v$ of $F$, the Langlands parameter
of the local component $\tilde{\Pi}(I(\theta))_v$ is the local
component $I(\theta)_v$. By the Jacquet-Langlands correspondence
\cite{JL}, this can be lifted to an automorphic representation,
denoted again by $\tilde{\Pi}(I(\theta))$ of $GL_1(D)(\A)$ provided
$\theta$ is admissible with respect to $D$.
 
\begin{definition}
An irreducible admissible representation $\Pi$ of $G(\A)$ is said to
be {\em dihedral}, if there exists a quadratic extension $E$ of $F$
and  an admissible idele class character $\theta$ of $E$ with respect
to $D$ such that  the local component $\Pi_v$  occurs in the
restriction of the local component $\tilde{\Pi}(I(\theta))_v$ to
$G(F_v)$. Further, at almost all finite places $v$ of $F$ where $D$ is
unramified,  the subgroup $SL(2, {\mathcal O}_{F_v})$ has a nonzero
fixed vector in the space of $\Pi_v$. We will say that $\Pi$ occurs in
the restriction of $\tilde{\Pi}(\theta)$ to $G(\A)$.
\end{definition}

Given a representation $\Pi$ of $G(\A)$,  let $m(\Pi)$ be the
multiplicity of $\Pi$ occuring in the space of cusp forms on $G(\A)$.

\begin{theorem}[Labesse-Langlands]\label{ll}
If $\Pi$ is a non-dihedral  cuspidal automorphic representation $\Pi$
  of $G(\A)$,  then for any $g\in \tilde{G}(\A)$, the representations
  $\Pi$ and   $\Pi^g$ occur with the same multiplicity in the space of
  cusp forms on $G(\A)$.
\end{theorem}

\subsubsection{Multiplicity formula} \label{ssection:mult}
We first define an auxiliary integer that occurs in  the multiplicity
 formula of Labesse and Langlands for the multiplicity $m(\Pi)$ of a
 representation of dihedral type in the space of cusp forms of
 $G(\A)$.

Given two parameters, $\tilde{\phi}, ~\tilde{\phi}':W_F\to GL_2(\C)$
call them to be weakly globally equivalent if for any place $v$ of
$F$, there exists a character $\chi_v$ of $F_v*$ such that
$\tilde{\phi}_v=\chi_v\tilde{\phi}_v'$; and globally equivalent if
there exists an idele class character $\chi$ of $C_F$ such that
$\tilde{\phi}=\chi\tilde{\phi}'$. Define the integer $d(\tilde{\phi})$
to be the number of weak global equivalence classes modulo global
equivalence.

Let $\Pi$ be a representation of $G(\A)$ occuring  in the restriction
of $\tilde{\Pi}(I(\theta))$ to $G(\A)$. Let $\tilde{\phi}=I(\theta)$
be the Langlands parameter associated to an admissible character
$\theta$ as defined above.  Let $\phi$ denote the projection of
$\tilde{\phi}$ to $PGL_2(\C)$. Define $d(\Pi)=d(I(\theta))$.  It
follows from the definition of $L$-indistinguishability, that for
$g\in \tilde{G}(\A)$,
\begin{equation}\label{inv:dpi}
d(\Pi)=d(\Pi^g).
\end{equation}

For each place $v$ of $F$, let $\phi_v: W_{F_v}\to PGL_2(\C)$  denote
the local component of the parameter $\phi$.   We have natural maps
$S_{\phi}\to S_{\phi_v}$, sending $S_{\phi}^0$ to
$S_{\phi_v}^0$. Given an element $s\in S_{\phi}$, we denote by $s_v$
both it's image in $S_{\phi_v}$ and in the group of connected
components $C_{\phi_v}$ of $S_{\phi_v}$.  We have,
\begin{theorem}[Multiplicity formula] \label{thm:mult}
Fix an additive character $\psi$ of $\A_F$.  Let  $\Pi$ be a
representation of $G(\A)$ occuring  in the restriction of
$\tilde{\Pi}(I(\theta))$ to $G(\A)$. With notation as above,  the
multiplicity of $\Pi$ in the space of automorphic forms of $G(\A)$ is
given by the formula:
\[ m(\Pi)=\frac{d(\Pi)}{[S_{\phi}/S_{\phi}^0]}\sum_{s\in
  S_{\phi}/S_{\phi}^0}\prod_{v\in {\Sigma_F}}<s_v, \Pi_v>, \] where
the local pairings  are normalized with respect to the additive
character $\psi_v$ of $F_v$, the local component of $\psi$ at $v$.

\end{theorem} 

We now make this formula more explicit.
\begin{definition} \cite[pages 42-44]{LL}
Let $\theta$ be an idele class character of $J_E/E^*$ admissible with
  respect to $D$.  The character $\theta$ is said to be of $type ~(b)$
  if $\theta/\theta^{\sigma}$ is a (non-trivial) quadratic
  character. In this case $C_{\phi}\simeq \Z/2\Z\times \Z/2\Z$.

If  $\theta$  is not of $type ~(b)$, then it is said to be of $type
  ~(a)$. Here  $C_{\phi}\simeq \Z/2\Z$, which we identify with
  $\{1,\epsilon\}$.

 Call $\theta$ to be of $type ~(a')$ if it is of $type~ (a)$,  and
  there exists a place $v$ of $F$ where $D$ does not split, such that
  $\theta_w/{\theta}^{\sigma}_w$ is a quadratic character of
  $E_w^*$. Here $w$ is the unique place of $E$ lying over the place
  $v$ of $F$.

The character $\theta$ is said to be of $type ~(a'')$ if it is of type
  (a) and not of type (a').
\end{definition}
 
Combining Theorem \ref{thm:mult} together with the local properties of
the pairing given by Proposition \ref{pairing:nondeg} and Equation
\ref{pairing:exceptional}, we obtain the following corollary
(Propositions 7.3 and 7.4 of \cite{LL}):
\begin{corollary}\label{cor:mult} 
Let $\Pi=\Pi(\theta)$ be a dihedral of $G(\A)$, for some admissible
  character $\theta$ of a quadratic extension $E$ of $F$ with respect
  to $D$.
\begin{enumerate} 
\item If $\theta$ is of $type ~(a')$, then
\[ m(\Pi) = \frac{d(\Pi)}{2}\prod_v<1,\Pi_v>.\]
\item  If $\theta$ is of $type ~(a'')$,  then
\[ m(\Pi) = \frac{d(\Pi)}{2}\left(\prod_v<1,\Pi_v>+\prod_v<\epsilon,
 \Pi_v>\right).\]
\item If $\theta$ is of $type ~(b)$, then
\[ m(\Pi) = \frac{d(\Pi)}{4}\prod_v<1,\Pi_v>.\] 
\end{enumerate}
\end{corollary}

\section{Proof of Theorem \ref{adelicconjtheorem}} \label{proof}
In order to prove the theorem, we  have to show  that
\[ L^2(\Gamma_K\backslash G_{\infty})\simeq  L^2(\Gamma_K^x\backslash
G_{\infty}), \] as $G_{\infty}$-modules.  Fix an unitary
representation $\pi$ of $G_{\infty}$. We  need to show  the equality
of multiplicities:
\begin{equation}\label{eqofmult}
m(\pi, \Gamma_K)=m(\pi, \Gamma_{K^x}).
\end{equation} 

Let ${\mathcal A}$ denote the equivalence classes of automorphic
representations of $G(\A)$. Since $G$ is simply connected and
semisimple, and $G_{\infty}$ is non-compact, by strong approximation
we obtain,
\[ L^2(\Gamma_K\backslash G_{\infty})\simeq \oplus_{\Pi\in {\mathcal
A}}m(\Pi)\text{dim}(\Pi_f^K)\Pi_{\infty}.\] Here $\Pi_f^K$ denotes the
space of $K$-invariants in the representation space for $\Pi_f$, and
$m(\Pi)$ is the multiplicity with which $\Pi$ occurs in the space of
automorphic forms of $G(\A)$.  Consequently,
\[ m(\pi, \Gamma_K)=\sum_{\{\Pi\in {\mathcal A}~\mid
  ~\Pi_{\infty}=\pi\}}m(\Pi)\text{dim}(\Pi_f^K).\] In order to
establish Equation \ref{eqofmult}, it is enough to show the following:
\begin{equation}  \label{eqn:glmult}
\sum_{\{\Pi\in {\mathcal A}~\mid
~\Pi_{\infty}=\pi\}}m(\Pi)\text{dim}(\Pi_f^K)= \sum_{\{\Pi\in
{\mathcal A}~\mid~ \Pi_{\infty}=\pi\}}m(\Pi)\text{dim}(\Pi_f^{K^x}).
\end{equation}

Fix a representation $\pi_{\infty}$ of $G_{\infty}$. Given a $\Pi$
with archimedean component $\Pi_{\infty}=\pi$, we will produce another
representation $\Pi'$ in the same $L$-packet satisfying,
\begin{itemize}
\item The archimedean component $\Pi'_{\infty}=\pi$.
\item The map $\Pi\mapsto \Pi'$ is injective, and

\begin{equation}\label{eqn:glmult2}
 m(\Pi)\text{dim}(\Pi_f^K)=m(\Pi')\text{dim}({\Pi_f'}^{K^x}).
\end{equation}

\end{itemize}

  By Equation \ref{eqn:glmult} it follows that $ m(\pi, \Gamma_K)\leq
m(\pi, \Gamma_{K^x})$. Reversing the argument we obtain Equation
\ref{eqn:glmult}.

If $\Pi$ is not dihedral,  let $\Pi'=\Pi^x$. By Theorem \ref{ll},
$m(\Pi)=m(\Pi^x)$. Clearly
$\text{dim}(\Pi_f^K)=\text{dim}({\Pi_f^x}^{K^x})$  and Equation
\ref{eqn:glmult2} follows.

Now suppose $\Pi$ is a dihedral associated to an idele class character
$\theta$ of a quadratic extension $E$ over $F$ of $type ~(a')$ or
$type ~(b)$.  Let  $\Pi'=\Pi^x$. By  parts (1) and (3) of Corollary
\ref{cor:mult} and Equation \ref{inv:dpi},
\[m(\Pi)=m(\Pi^x),\]
and  Equation \ref{eqn:glmult2} follows.

Hence we have to only consider dihedral representations $\Pi$ of the
form $\pi(\theta)$ (with archimedean component  $\pi$), where $\theta$
is an idele class character of a quadratic extension $E/F$ of $type
~(a'')$.  Let $w_0$ be the unique place of $E$ lying over the
distinguished place $v_0$ of $F$. By \ref{pairing:crucial}, at the
place $v_0$, there are two elements in the L-packet
\[L(I(\theta_{w_0}))=\{\pi(I(\theta_{w_0}))^+, \pi(I(\theta_{w_0}))^-\},\] 
which pairs perfectly with the group of connected components
$C_{I(\theta_{w_0})}\simeq \Z/2\Z=\{1,\epsilon\}$, as given by
Equation \ref{pairing:crucial}.

Let $T$ be the finite subset of finite places $v$ of $F$,  such that
\[<\epsilon, \pi_v>\neq <\epsilon, \pi_v^{x_v}>,\] where $x_v$ is the
component at $v$ of the adele $x\in \tilde{G}(\A_f)$.

Define the local components of $\Pi'$ at any place $v$ of $F$ as,
\[\Pi_v'=\begin{cases} \Pi_v^x \quad \text{if $v\neq v_0$},\\
\pi(I(\theta_{w_0}))^+ \quad \text{if $v=v_0$ and $|T|$ is even},\\
\pi(I(\theta_{w_0}))^-  \quad \text{if $v=v_0$ and $|T|$ is odd}.
\end{cases}
\]
Since the representations $\pi(\theta)_{v_0}^{\pm}$ are conjugate
under $D^*_{v_0}$,  both $\Pi$ and $\Pi'$ are in the same $L$-packet,
and  by Equation \ref{inv:dpi},
\[d(\Pi)=d(\Pi').\]
It follows from part (2) of Corollary \ref{cor:mult}, that
\[m(\Pi)=m(\Pi').\] 
By Hypothesis {\bf H3}, $K_{v_0}$ is a normal subgroup of $D^*_{v_0}$.
Hence,
\[ \text{dim}(\Pi_f^K)=\text{dim}({\Pi_f'}^{K^x}).\]
From these  equations, we conclude that Equation
\ref{eqn:glmult2} is satisfied. 

This concludes the proof of Theorem \ref{adelicconjtheorem}.

\end{document}